\documentclass{article}

\usepackage{latexsym,amssymb}
\usepackage[noadjust]{cite}
\input amssym.def \input amssym

\newtheorem{te}{Theorem}[section]

\newtheorem{de}[te]{Definition}
\newtheorem{lm}[te]{Lemma}

\newtheorem{pp}[te]{Proposition}
\newtheorem{co}[te]{Corollary}
\newtheorem{ex}[te]{Example}

\def\dokaz{\noindent{\bf Proof. }}
\def\kraj{\hfill $\Box$ \par \vspace*{2mm} }
\def\widemid{\hspace{1mm}\widetilde{\mid}\hspace{1mm}}
\def\nwidemid{\hspace{1mm}\widetilde{\nmid}\hspace{1mm}}

\newcommand{\zve}[1]{{{}^*\hspace{-0.5mm}#1}}
\newcommand{\zvez}[1]{{{}^*\hspace{-1mm}#1}}
\def\zvepar{\hspace{1mm}\zvez\parallel\hspace{1mm}}
\def\zvemid{\hspace{1mm}\zvez\mid\hspace{1mm}}
\def\nzvemid{\hspace{1mm}\zvez\nmid\hspace{1mm}}

\begin{document}

\begin{center}
           {\huge \bf Divisibility classes of ultrafilters and their patterns}\\[2mm]
{\bf Boris  \v Sobot}\\[2mm]
{\small  Department of Mathematics and Informatics, University of Novi Sad,\\
Trg Dositeja Obradovi\'ca 4, 21000 Novi Sad, Serbia\\
e-mail: sobot@dmi.uns.ac.rs, ORCID: 0000-0002-4848-0678}
\end{center}

\begin{abstract}
A divisibility relation on ultrafilters on the set $\mathbb{N}$ of natural numbers is defined as follows: ${\cal F}\hspace{1mm}\widetilde{\mid}\hspace{1mm}{\cal G}$ if and only if every set in $\cal F$ upward closed for divisibility also belongs to $\cal G$. Previously we isolated basic classes: powers of prime ultrafilters, and described the pattern of an ultrafilter, measuring the quantity of members of each basic class dividing a given ultrafilter. In this paper we define a topology on the set of basic classes which will allow us to calculate the pattern of the limit of a $\widetilde{\mid}$-increasing chain of ultrafilters. Using this we characterize which patterns can actually appear as patterns of an ultrafilter. Defining the $=_\sim$-divisibility classes by identifying mutually divisible ultrafilters, in the respective quotient order $(\beta\mathbb{N}/=_\sim,\widemid)$ we identify singleton classes and consider their patterns. Finally, we give a sufficient condition for a $=_\sim$-divisibility class to have an immediate predecessor.
\end{abstract}



Ultrafilters are often used in place of choice principles in various constructions. For this purpose it is sometimes useful to have an ultrafilter with special properties. A typical example is the ultrafilter proof of Hindman's theorem, which uses an idempotent ultrafilter (see \cite{HS} for this and many other examples). Hence, for constructions involving number-theoretic concepts it may be useful to have ultrafilters on the set $\mathbb{N}$ of natural numbers with some properties derived from divisibility. This was the motivation for defining a divisibility relation on ultrafilters in \cite{So1}.

In this paper we continue the structural analysis of ultrafilters from the point of view of their place within the divisibility quasiorder, initiated in \cite{So8}. This is accomplished via so-called patterns: they reflect the structure of nonstandard generators of an ultrafilter.

The first of the two main results of the paper is Theorem \ref{limpat}, calculating the pattern of the limit of a sequence of ultrafilters. It provides a way of controlling the pattern within recursive constructions. An example of such a construction leads to the other main result, Corollary \ref{characterization}, which characterizes the patterns which can be realized by an ultrafilters. This, in turn, provides information on the structure of nonstandard generators of an ultrafilter. Actually, what we obtain is a general result on the number of powers of primes generating a given basic class that can be exact factors of a nonstandard number (Corollary \ref{nonst}).

\section{Introduction: divisibility and patterns}

Whenever we write ``ultrafilter" it will be understood, unless otherwise stated, that it is an ultrafilter on the set of natural numbers $\mathbb{N}$. Let $\beta\mathbb{N}$ be the set of all ultrafilters on $\mathbb{N}$. For each $n\in\mathbb{N}$, the principal ultrafilter $\{A\subseteq\mathbb{N}:n\in A\}$ is identified with $n$. Considering the topology with base sets $\overline{A}=\{{\cal F}\in\beta\mathbb{N}:A\in{\cal F}\}$, $\beta\mathbb{N}$ is an extension of the discrete space $\mathbb{N}$, called the Stone-\v Cech compactification of $\mathbb{N}$. 

Every function $f:\mathbb{N}\rightarrow\mathbb{N}$ can be uniquely extended to a continuous $\widetilde{f}:\beta\mathbb{N}\rightarrow\beta\mathbb{N}$. For ${\cal F}\in\beta\mathbb{N}$, $\widetilde{f}({\cal F})$ is the ultrafilter generated by $\{f[A]:A\in{\cal F}\}$. 
In \cite{So1} a similar way to extend the divisibility relation $\mid$ to $\beta\mathbb{N}$ was proposed: if
$$A\hspace{-1mm}\uparrow=\{n\in {\mathbb{N}}:(\exists a\in A)a\mid n\}$$
for $A\subseteq\mathbb{N}$, let ${\cal F}\hspace{1mm}\widetilde{\mid}\hspace{1mm}{\cal G}$ if for every $A\in{\cal F}$ holds $A\hspace{-1mm}\uparrow\in{\cal G}$. A more practical equivalent condition is the following:
$${\cal F}\hspace{1mm}\widetilde{\mid}\hspace{1mm}{\cal G}\Leftrightarrow{\cal F}\cap{\cal U}\subseteq{\cal G},$$
where
$${\cal U}=\{A\in P({\mathbb{N}})\setminus\{\emptyset\}:A=A\hspace{-1mm}\uparrow\}$$
is the family of nonempty upwards-closed sets. The relation $\widemid$ is a quasiorder, and defining
$${\cal F}=_\sim{\cal G}\Leftrightarrow{\cal F}\widemid{\cal G}\land{\cal G}\widemid{\cal F},$$
$\widemid$ can be thought of as a partial order on equivalence classes from $\beta\mathbb{N}/\hspace{-0.5mm}=_\sim$, which we denote simply by $[{\cal F}]$. 
Various aspects of $(\beta\mathbb{N}/\hspace{-1mm}=_\sim,\widemid)$ were considered in papers \cite{So1}--\cite{So9}.

Parallel to $\beta\mathbb{N}$ we can observe a nonstandard extension of $\mathbb{N}$, denoted by $\zve{\mathbb{N}}$. All the prerequisites about nonstandard extensions that we need can be found in \cite{G2}, \cite{DGL} or in the previous papers on the subject \cite{So5}, \cite{So8} and \cite{So9}. Here we just emphasize that, in a nonstandard model, every subset of $\mathbb{N}$, every function and every relation on $\mathbb{N}$ has its star-version on $\zve{\mathbb{N}}$, satisfying the same formulas with bounded quantifiers; this is called the Transfer principle. For example, the divisibility relation $\mid$ on $\mathbb{N}$ is extended to a relation $\zvemid$ on $\zve{\mathbb{N}}$.

The connection between nonstandard and ultrafilter extensions of $\mathbb{N}$ is obtained as follows: for every $x\in\zve{\mathbb{N}}$ we identify the 1-type of $x$ over $\mathbb{N}$ with the corresponding ultrafilter:
$$tp(x/\mathbb{N}):=\{A\subseteq\mathbb{N}:x\in\zve A\}.$$
The element $x$ is then called a generator of that ultrafilter ${\cal F}:=\{A\subseteq\mathbb{N}:x\in\zve A\}$ and we write $x\models{\cal F}$. This agrees with extensions of functions $f:\mathbb{N}\rightarrow\mathbb{N}$: $tp(\zve f(x)/\mathbb{N})=\widetilde{f}(tp(x/\mathbb{N}))$ for every $x\in\zve{\mathbb{N}}$.

Thus, $\zve{\mathbb{N}}$ serves as a mirror to $\beta\mathbb{N}$: properties of generators of an ultrafilter $\cal F$ reflect properties of $\cal F$. Additional assumptions on $\zve{\mathbb{N}}$ give us a sharper image, simplifying some definitions and proofs. An example is the following basic result, providing a way to translate statements about divisibility from $\beta\mathbb{N}$ to $\zve{\mathbb{N}}$ and vice versa.

\begin{pp}[\cite{So5}, Theorem 3.4]\label{ekviv}
If $\zve{\mathbb{N}}$ is ${\goth c}^+$-saturated, the following conditions are equivalent for every two ultrafilters ${\cal F},{\cal G}\in\beta \mathbb{N}$:

(i) ${\cal F}\widemid{\cal G}$;

(ii) there are generators $x\models{\cal F}$ and $y\models{\cal G}$ such that $x\zvemid y$;

(iii) for every generator $x\models{\cal F}$ there is a generator $y\models{\cal G}$ such that $x\zvemid y$;

(iv) for every generator $y\models{\cal G}$ there is a generator $x\models{\cal F}$ such that $x\zvemid y$.
\end{pp}

{\bf Throughout this text we assume that we work with a nonstandard extension which is ${\goth c}^+$-saturated and satisfies Henson's Isomorphism property.} This property was defined in \cite{He} and further studied in \cite{DN3}. Such extensions exist in ZFC, and in them all the sets $[1,z]:=\{a\in\zve{\mathbb{N}}:a\leq z\}$ for $z\in\zve{\mathbb{N}}\setminus\mathbb{N}$ have the same cardinality that we denote $\infty$. Using ${\goth c}^+$-saturation it is easy to see that every $z\in\zve{\mathbb{N}}\setminus\mathbb{N}$ is preceded by at least one generator of each ultrafilter, so $\infty\geq 2^{\goth c}$.

The effect of this assumption will be that we will deal with only one infinite cardinality of practically all the sets we work with. But the sole fact that such an extension exists shows that possible different cardinalities (in other nonstandard extensions) do not reflect properties of ultrafilters, but only peculiarities of the extension itself. Hence restricting ourselves in this way will not narrow down our view on ultrafilters.

We denote $\mathbb{N}_\infty=\omega\cup\{\infty\}$. Topologically, $\mathbb{N}_\infty$ is the one-point compactification of the discrete space on $\omega$.

We recall some more notation and previous results. By $\mathbb{P}$ we denote the set of prime numbers. Ultrafilters containing $\mathbb{P}$ are called prime; they are divisible only by 1 and themselves and their generators are exactly the prime numbers from $\zve{\mathbb{P}}$. For $A,B\subseteq\mathbb{N}$ denote $A^c=\mathbb{N}\setminus A$ and $A\cdot B=\{ab:a\in A\land b\in B\land gcd(a,b)=1\}$. For $n\in\mathbb{N}$, let $A^{(n)}=\underbrace{A\cdot A\cdot\dots\cdot A}_n$. 

Now denote $\mathbb{P}^{exp}:=\{p^x:p\in\mathbb{P}\land x\in\mathbb{N}\}$. Elements of $\zve{\mathbb{P}^{exp}}$ generate some ultrafilters; when we identify those of them that are $=_\sim$-equivalent, we get so-called basic classes, denoted by ${\cal P}^u$. They serve as powers of prime ultrafilters ${\cal P}\in\mathbb{P}$. 

\begin{de}
Let ${\cal P}\in\overline{\mathbb{P}}$. The relation $\approx_{\cal P}$ on $\{(p,a):p\models{\cal P}\land a\in\zve{\mathbb{N}}\}$ is defined as follows:
$$(p,a)\approx_{\cal P}(q,b)\mbox{ if and only if }tp(p^a/\mathbb{N})=_\sim tp(q^b/\mathbb{N}).$$
Then $\approx_{\cal P}$ is an equivalence relation; let
$${\cal E}_{\cal P}=\{[(p,a)]_{\approx_{\cal P}}:p\models{\cal P}\land a\in\zve{\mathbb{N}}\}$$
be the set of its equivalence classes.

On ${\cal E}_{\cal P}$ we define the relation:
\begin{eqnarray*}
u\prec_{\cal P}v &\mbox{ if and only if } & u\neq v\mbox{ and for some } p\models{\cal P}\mbox{ and some }a,b\in\zve{\mathbb{N}}\\
&& \mbox{ holds }(p,a)\in u,(p,b)\in v\mbox{ and }a<b.
\end{eqnarray*}
We write $u\preceq_{\cal P}v$ if $u\prec_{\cal P}v$ or $u=v$.

Families of ultrafilters of the form ${\cal P}^u:=\{tp(p^a/\mathbb{N}):(p,a)\in u\}$ for some ${\cal P}\in\overline{\mathbb{P}}$ and $u\in{\cal E}_{\cal P}$ are called basic classes. By ${\cal B}$ we denote the set of all basic classes.
\end{de}

Of course, $({\cal E}_{\cal P},\prec_{\cal P})$ is a linear order. For $p\in\mathbb{P}$ there is only one ``infinite" element of ${\cal E}_p$, and the corresponding basic class is denoted by $p^\omega$. When we write $A\in{\cal P}^u$ for some $A\subseteq\mathbb{N}$, it means that $A\in{\cal F}$ for all ${\cal F}\in{\cal P}^u$. Generators $p^a$ of some ${\cal F}\in{\cal P}^u$ are also called generators of ${\cal P}^u$ and we write $p^a\models{\cal P}^u$. 

If ${\cal P}\in\overline{\mathbb{P}}\setminus\mathbb{P}$ we distinguish ${\cal P}^u\in{\cal B}$ of the first kind (such that, for every $p\models{\cal P}$, the set $u_p:=\{a\in\zve{\mathbb{N}}:p^a\models{\cal P}^u\}$ is a singleton) and those of the second kind (such that, for every $p\models{\cal P}$, $u_p$ is a union of consecutive galaxies - blocks of order type $(\mathbb{Z},<)$ in $\zve{\mathbb{N}}$). The generators of basic classes of the first kind are of the form $p^{\zve f(p)}$ for some $f:\mathbb{N}\rightarrow\mathbb{N}$, so such classes have immediate predecessors and immediate successors (generated by $p^{\zve f(p)-1}$ and $p^{\zve f(p)+1}$ respectively). Basic classes ${\cal P}^u$ of the second kind have neither, and we can represent $u$ as a supremum of an increasing sequence, or an infimum of a decreasing sequence in $({\cal E}_{\cal P},\prec_{\cal P})$. For more information on $({\cal E}_{\cal P},\prec_{\cal P})$ see \cite{So8}.

Let $\langle p_n:n\in\zve{\mathbb{N}}\rangle$ be the increasing enumeration of $\zve{\mathbb{P}}$. In any nonstandard extension a generalization of the Fundamental theorem of arithmetic holds, saying that every $x\in\zve{\mathbb{N}}$ can be written uniquely as $\prod_{n\leq z}p_n^{h(n)}$ for some $z\in\zve{\mathbb{N}}$ and some internal function $h:[1,z]\rightarrow\zve{\mathbb{N}}$ so that $h(z)>0$. Also, $p^a\zvepar x$ means that $x=p^ax'$ where $p\nzvemid x'$; we say that $p^a$ is an exact divisor of $x$.

\begin{de}\label{defpattern}
Let ${\cal A}$ be the set of all functions $\alpha:{\cal B}\rightarrow\mathbb{N}_\infty$ such that $\sum_{k\in\omega+1}\alpha(p^k)\leq 1$ for $p\in\mathbb{P}$. Elements $\alpha\in{\cal A}$ are called patterns.


For any $x=\prod_{n\leq z}p_n^{h(n)}\in\zve{\mathbb{N}}$ define its pattern $\alpha_x\in{\cal A}$ as follows. For each basic ${\cal P}^u\in{\cal B}$, let $\alpha_x({\cal P}^u):=|\{(p,a):p^a\models{\cal P}^u\land p^{a}\zvepar x\}|$.
\end{de}

Thus, $\alpha_x({\cal P}^u)$ counts exact divisors of $x$ which are generators of ${\cal P}^u$. It was proved in \cite{So8}, Theorem 3.3 that this cardinality is indeed always in $\mathbb{N}_\infty$. The condition $\sum_{k\in\omega+1}\alpha(p^k)\leq 1$ for $p\in\mathbb{P}$ reflects the fact that only one power of $p$ can be an exact divisor of any $x$.

We sometimes write $\alpha=\{({\cal P}_i^{u_i},k_i):i\in I\}$, meaning that $\alpha({\cal P}_i^{u_i})=k_i$ and $\alpha({\cal Q}^v)=0$ for basic classes ${\cal Q}^v$ distinct from all ${\cal P}_i^{u_i} $.

Finally, for ${\cal F}\in\beta\mathbb{N}$, $\alpha_{\cal F}:=\alpha_x$ for any generator $x$ of $\cal F$. It was proved in \cite{So8}, Theorem 3.12 that this definition does not depend on the choice of $x$.

\section{Limits of $\widemid$-chains}

When we need an ultrafilter $\cal F$ with some prescribed divisibility-related properties, one natural way to construct it would be to approach it ``from below". In other words, we can try to approximate $\cal F$ by elements of a $\widemid$-increasing chain. We now formalize this idea.

If $\langle{\cal F}_i:i\in I\rangle$ is a sequence of ultrafilters and $\cal H$ is an ultrafilter on $I$, the limit of the sequence is the ultrafilter defined by: $A\in\lim_{i\rightarrow\mathcal{H}}\mathcal{F}_i$ if and only if $\{i\in I:A\in\mathcal{F}_i\}\in\mathcal{H}$. 

\begin{pp}[\cite{So5}, Lemma 4.1.]\label{sublimit}
Let $\langle[\mathcal{F}_i]:i\in I\rangle$ be a chain in the order $(\beta \mathbb{N}/\hspace{-1mm}=_\sim,\widemid)$.

(a) It has the least upper bound $[\mathcal{G}_U]$, obtained as $\mathcal{G}_U=\lim_{i\rightarrow\mathcal{H}}\mathcal{F}_i$ for any ultrafilter $\mathcal{H}$ containing all final segments of $I$.

(b) It has the greatest lower bound $[\mathcal{G}_L]$, obtained as $\mathcal{G}_L=\lim_{i\rightarrow\mathcal{H}}\mathcal{F}_i$ for any ultrafilter $\mathcal{H}$ containing all initial segments of $I$.

(c) $\mathcal{G}_U\cap\mathcal{U}=\bigcup_{i\in I}(\mathcal{F}_i\cap\mathcal{U})$ and $\mathcal{G}_L\cap\mathcal{U}=\bigcap_{i\in I}(\mathcal{F}_i\cap\mathcal{U})$.
\end{pp}

So in both cases the $=_\sim$-equivalence class of the limit does not depend on the ultrafilter ${\cal H}$ we take the limit over (as long as it contains all initial/final segments of $I$). Hence, if $\langle[\mathcal{F}_\xi]:\xi\in\gamma\rangle$ is a well-ordered $\widemid$-increasing chain, we can write $[\mathcal{G}_U]:=\lim_{\xi\rightarrow\gamma}[\mathcal{F}_\xi]$, meaning that we take $\mathcal{G}_U$ to be any ultrafilter from the $=_\sim$-equivalence class containing $\lim_{\xi\rightarrow{\cal H}}\mathcal{F}_\xi$ for all $\cal H$ containing all final segments of $\gamma$.

Applying the proposition above to the linear order of basic classes yields the following corollary.

\begin{pp}[\cite{So8}, Lemma 2.7]\label{supinf}
For every ${\cal P}\in\overline{\mathbb{P}}$, every increasing sequence in $({\cal E}_{\cal P},\prec_{\cal P})$ has a supremum and every decreasing sequence has an infimum.
\end{pp}

In particular there are the greatest element and the smallest infinite element in ${\cal E}_{\cal P}$; the corresponding basic classes are denoted by ${\cal P}^{max}$ and ${\cal P}^{\omega}$.

Lemma 4.4 from \cite{So5} states that every well-ordered $\widemid$-chain $\langle{\cal F}_\xi:\xi<\gamma\rangle$ ``reflects" in $\zve{\mathbb{N}}$: there is a $\zvemid$-chain of generators $x_\xi\models{\cal F}_\xi$. We adjust its proof to show the following.

\begin{lm}\label{midniz}
Let $\gamma<{\goth c}^+$, $\langle x_\xi:\xi<\gamma\rangle$ is a $\zvemid$-chain, and let ${\cal G}\in\beta\mathbb{N}$ be divisible by ${\cal F}_\xi:=tp(x_\xi/\mathbb{N})$ for all $\xi<\gamma$.

(a) There is $y\models{\cal G}$ such that $x_\xi\zvemid y$ for all $\xi<\gamma$.

(b) For any $y\models{\cal G}$ as in (a), there is $z\in\zve{\mathbb{N}}$ such that $x_\xi\zvemid z$ for all $\xi<\gamma$, $z\zvemid y$ and $[tp(z/\mathbb{N})]$ is the least upper bound of $\{[{\cal F}_\xi]:\xi<\gamma\}$.
\end{lm}

\dokaz (b) Let $F_0:=\{A:A^c\in{\cal U}\land(\forall\xi<\gamma)x_\xi\in\zve A\}$ and $F:=\{\zve A:A\in F_0\}\cup\{\{z\in\mathbb{N}:x_\xi\zvemid z\}:\xi<\gamma\}\cup\{\{z\in\mathbb{N}:z\zvemid y\}\}$. We are looking for an element in $\bigcap F$. Since $F$ is a family of internal sets, by ${\goth c}^+$-saturation it suffices to show that it has the finite intersection property. So let $x_{\xi_1},x_{\xi_2},\dots,x_{\xi_k}$ and $A_1,A_2,\dots,A_l\in F_0$ be given. Let $\xi=\max\{\xi_1,\xi_2,\dots,\xi_k\}$. Then $x_\xi$ itself is the element divisible by all $x_{\xi_i}$, dividing $y$ and contained in all $\zve A_j$.

(a) is proved in a similar way.\kraj

\section{Patterns of limits}

\begin{de}
A set $A$ is convex if, for all $a,b\in A$, $a\mid c$ and $c\mid b$ imply $c\in A$. 

Let $\cal C$ denote the family of all nonempty convex subsets of $\mathbb{P}^{exp}$. 

Define $\overline{A}:=\{{\cal P}^u\in{\cal B}:A\in{\cal P}^u\}$ for $A\in{\cal C}$.
\end{de}

This duplicate notation ($\overline{A}$ as a set of basic classes and $\overline{A}=\{{\cal F}\in\beta\mathbb{N}:A\in{\cal F}\}$ as a set of ultrafilters) should not cause a confusion since, for $A\in{\cal C}$ and ${\cal P}^u\in{\cal B}$,
\begin{eqnarray*}
 & & {\cal P}^u\in\overline{A}\\
\mbox{if and only if } & & {\cal F}\in\overline{A}\mbox{ for every }{\cal F}\in{\cal P}^u\label{overlineA}\\
\mbox{if and only if } & & {\cal F}\in\overline{A}\mbox{ for some }{\cal F}\in{\cal P}^u.
\end{eqnarray*}
Of course, if $A\in{\cal C}$, then $\zve A\subseteq\zve{\mathbb{P}}^{exp}$.

\begin{lm}\label{topology}
(a) The sets $\overline{A}$ for $A\in{\cal C}$ constitute a basis for some topology $\cal O$ on $\cal B$.

(b) The space $({\cal B},{\cal O})$ is Hausdorff.

(c) The space $({\cal B},{\cal O})$ is compact.
\end{lm}

\dokaz (a) Obvious, since the intersection of convex sets is also convex.

(b) Let ${\cal P}^u,{\cal Q}^v\in{\cal B}$ be distinct. If ${\cal P}\neq{\cal Q}$, there is $A\subseteq\mathbb{P}$ such that $A\in{\cal P}\setminus{\cal Q}$, so $A^{exp}:=\{p^n:p\in A\land n\in\mathbb{N}\}\in{\cal P}^u\cap{\cal C}\setminus{\cal Q}^v$ and $(\mathbb{P}\setminus A)^{exp}\in{\cal Q}^v\cap{\cal C}\setminus{\cal P}^u$. Otherwise, if for example $u\prec_{\cal P}v$, there is some $A\in{\cal P}^v\cap{\cal U}\setminus{\cal P}^u$, so $B:=A\cap\mathbb{P}^{exp}\in{\cal P}^v\cap{\cal C}\setminus{\cal P}^u$ and $\mathbb{P}^{exp}\setminus B\in{\cal P}^u\cap{\cal C}\setminus{\cal P}^v$.

(c) The set $\overline{\mathbb{P}^{exp}}$ is closed in $\beta\mathbb{N}$ and thus compact. The map $\varphi:\overline{\mathbb{P}^{exp}}\rightarrow{\cal B}$, defined by $\varphi({\cal F})={\cal P}^u$ if ${\cal F}\in{\cal P}^u$, is continuous and onto so $({\cal B},{\cal O})$ is compact too.\kraj

In \cite{So8} a pattern $\alpha$ was called $\cal U$-closed if, whenever ${\cal P}^u\in{\cal B}$ and $n\in\mathbb{N}$, $\sum_{{\cal Q}^v\in\overline{A}}\alpha({\cal Q}^v)\geq n$ for every $A\in{\cal P}^u\cap{\cal U}$ implies $\sum_{w\succeq_{\cal P}u}\alpha({\cal P}^w)\geq n$. It was proven that the pattern $\alpha_{\cal F}$ of every ultrafilter is $\cal U$-closed. Now we define a strengthening of $\cal U$-closedness. 

\begin{de}
We will say that $\alpha\in{\cal A}$ is $\cal C$-closed if:

(C1) for every $A\in{\cal C}$, the value $\sum_{{\cal Q}^v\in\overline{A}}\alpha({\cal Q}^v)$ belongs to the set $\mathbb{N}_\infty$;


(C2) whenever ${\cal P}^u\in{\cal B}$ and $n\in\mathbb{N}$, $\sum_{{\cal Q}^v\in\overline{A}}\alpha({\cal Q}^v)\geq n$ for every $A\in{\cal P}^u\cap{\cal C}$ implies $\alpha({\cal P}^u)\geq n$.
\end{de}

Clearly, for $\alpha\in{\cal A}$ satisfying (C1), the condition (C2) is equivalent to any of the following two:

(C2') for all ${\cal P}^u\in{\cal B}$,
$$\alpha({\cal P}^u)=\min_{A\in{\cal P}^u\cap{\cal C}}\sum_{{\cal Q}^v\in\overline{A}}\alpha({\cal Q}^v);$$

(C2'') whenever $\alpha({\cal P}^u)<\infty$, there is some $A\in{\cal P}^u\cap{\cal C}$ such that $\alpha({\cal Q}^v)=0$ for all ${\cal Q}^v\in\overline{A}\setminus\{{\cal P}^u\}$.

We do not prove now that $\cal C$-closedness implies $\cal U$-closedness, since we will deduce this directly from further results as Corollary \ref{CimpliesU}.

Note also that (C2) does not affect basic classes $p^k$ for $p\in\mathbb{P}$ and $k\in\mathbb{N}$. Namely, such basic classes are singletons and they have a singleton neighborhood $\{p^k\}$.

\begin{lm}\label{cclosed}
For every ${\cal F}\in\beta{\mathbb{N}}$, the pattern $\alpha_{\cal F}$ is $\cal C$-closed.
\end{lm}

\dokaz Let $x\models{\cal F}$ be arbitrary, so $\alpha_{\cal F}=\alpha_x$. For $A\in{\cal C}$, the sum $\sum_{{\cal Q}^v\in\overline{A}}\alpha_x({\cal Q}^v)$ equals the cardinality of the set of exact divisors of $x$ which belong to $\zve A$. However, the set $Q=\{q\in\zve A:q\zvepar x\}$ is internal and bounded, so there is a bijection between $Q$ and the interval $[1,z]$ for some $z\in\zve{\mathbb{N}}$. Thus the cardinality of $Q$ is either finite or $\infty$, which proves (C1).

To show (C2), assume that $\sum_{{\cal Q}^v\in\overline{A}}\alpha_x({\cal Q}^v)\geq n$ for every $A\in{\cal P}^u\cap{\cal C}$. Denote $B_A=\{(q_1,q_2,\dots,q_n)\in\zve A^n:(\forall i\neq j)q_i\neq q_j\land(\forall i)q_i\zvepar x\}$. We prove that the family $F:=\{B_A:A\in{\cal P}^u\cap{\cal C}\}$ has the finite intersection property. Since this family is closed for finite intersections, we need only show that each $B_A$ is nonempty. The condition $\sum_{{\cal Q}^v\in\overline{A}}\alpha_x({\cal Q}^v)\geq n$ implies that there are distinct $z_i\models{\cal Q}_i^{v_i}$ for $1\leq i\leq n$ and some ${\cal Q}_i^{v_i}\in\overline{A}$ such that $z_i\zvepar x$. Thus $(z_1,z_2,\dots,z_n)\in B_A$, $F$ has the finite intersection property, so by ${\goth c}^+$-saturation we obtain $(q_1,q_2,\dots,q_n)\in\bigcap F$. The elements $q_i\in\zve{\mathbb{P}^{exp}}$ belong to $\zve A$ for all $A\in{\cal P}^u\cap{\cal C}$, so by Lemma \ref{topology}(b) they are generators of ${\cal P}^u$. Since they are are distinct and $q_i\zvepar x$ for $1\leq i\leq n$, it follows that $\alpha_x({\cal P}^u)\geq n$.\kraj

\begin{ex}\label{UCclosed}
Let us show that $\cal U$-closedness, even in conjuction with the condition (C1), does not imply (C2). Fix some ${\cal Q}\in\overline{\mathbb{P}}\setminus\mathbb{P}$ and take $\alpha=\{({\cal P},\infty):{\cal P}\in\overline{\mathbb{P}}\setminus(\mathbb{P}\cup\{{\cal Q}\})\cup\{({\cal P}^{max},\infty):{\cal P}\in\overline{\mathbb{P}}\setminus\mathbb{P}\}$. Every infinite $A\subseteq\mathbb{P}$ such that $A\in{\cal Q}$ belongs to ${\cal C}$, $\sum_{{\cal P}\in\overline{A}}\alpha({\cal P})=\infty$, but $\alpha({\cal Q})=0$, so $\alpha$ does not satisfy (C2). On the other hand, it satisfies (C1): each neighborhood $\overline{A}$ either contains no ${\cal P}^u$ for which $\alpha({\cal P}^u)>0$, or we have $\alpha({\cal P}^u)=\infty$ for any such ${\cal P}^u$. Also, $\alpha$ is $\cal U$-closed since: (a) for every ${\cal P}\in\overline{\mathbb{P}}\setminus\mathbb{P}$ and every $u\in{\cal E}_{\cal P}$, $\sum_{w\succeq_{\cal P} u}\alpha({\cal P}^u)\geq\alpha({\cal P}^{max})=\infty$, and (b) for $p\in\mathbb{P}$ and $k\in\mathbb{N}$, there are neighborhoods $A_k:=\{p^m:m\geq k\}\in{\cal C}$ of $p^k$ and $p^\omega$ such that ${\cal P}^u\notin\overline{A_k}$ for all ${\cal P}\neq p$.
\end{ex}

Let $\langle x_i:i\in I\rangle$ be a sequence in a topological space $X$, $y\in X$ and let $\cal H$ be an ultrafilter on $I$. Then ${\cal H}$-$\lim x_i=y$ means that, for every neighborhood $D$ of $y$, $\{i\in I:x_i\in D\}\in{\cal H}$. In particular, if $X=\beta\mathbb{N}$, ${\cal H}$-$\lim{\cal G}_i$ is the same as $\lim_{i\rightarrow{\cal H}}{\cal G}_i$ used above.

\begin{lm}\label{limult}
For every sequence $\langle x_i:i\in I\rangle$ of elements of $\mathbb{N}_\infty$ and every ultrafilter $\cal H$ on $I$,
$${\cal H}\mbox{-}\lim x_i = \min_{S\in{\cal H}}\sup_{i\in S}x_i = \sup_{S\in{\cal H}}\min_{i\in S}x_i.$$
\end{lm}

\dokaz We consider two cases. The first one: there is $n\in\mathbb{N}$ such that $\{i\in I:x_i=n\}\in{\cal H}$. Clearly, in this case all the given expressions equal $n$.

Otherwise, for every $S\in{\cal H}$ we have $\sup_{i\in S}x_i=\infty$, so $\min_{S\in{\cal H}}\sup_{i\in S}x_i=\infty$. Furthermore, for any $n\in\mathbb{N}$, $\{i\in S:x_i>n\}\in{\cal H}$, so $\sup_{S\in{\cal H}}\min_{i\in S}x_i=\infty$ as well. Also, neighborhoods of $\infty$ contain sets of the form $\{m:m>n\}$, so ${\cal H}$-$\lim x_i=\infty$.\kraj

We already mentioned that, for a given $\widemid$-chain $\langle{\cal G}_i:i\in I\rangle$, the $=_\sim$-equivalence class of ${\cal F}:=\lim_{i\rightarrow{\cal H}}{\cal G}_i$ does not depend on the choice of a particular ${\cal H}$. However, in order to control the pattern of the particular ultrafilter ${\cal F}$ that we obtain, we need to know what ${\cal H}$ is exactly. The result we obtain is true even if $\langle{\cal G}_i:i\in I\rangle$ is not a $\widemid$-chain; all we need is that $\lim_{i\rightarrow{\cal H}}{\cal G}_i$ exists.

\begin{te}\label{limpat}
Let $\{{\cal G}_i:i\in I\}$ be a family of ultrafilters on $\mathbb{N}$, let $\cal H$ be an ultrafilter on $I$, and ${\cal F}=\lim_{i\rightarrow{\cal H}}{\cal G}_i$. For any ${\cal P}^u\in{\cal B}$:
\begin{equation}\label{eqlimpat}
\alpha_{\cal F}({\cal P}^u)=\min_{A\in{\cal P}^u\cap{\cal C}}{\cal H}\mbox{-}\lim\sum_{{\cal Q}^v\in\overline{A}}\alpha_{{\cal G}_i}({\cal Q}^v).
\end{equation}
\end{te}

\dokaz Let us first prove the $\leq$ inequality. Assume the opposite; by Lemma \ref{limult} there are $A\in{\cal P}^u\cap{\cal C}$, $S\in{\cal H}$ and some $l\in\mathbb{N}$ such that for every $i\in S$
$$\alpha_{\cal F}({\cal P}^u)\geq l>\sum_{{\cal Q}^v\in\overline{A}}\alpha_{{\cal G}_i}({\cal Q}^v).$$
Let
$$B_A:=\bigcup_{k<l}(A^{(k)}\cdot\{n\in\mathbb{N}:\neg(\exists a\in A)a\parallel n\})$$
be the set consisting of all elements of $\mathbb{N}$ with less than $l$ exact divisors from $A$. Then $B_A\in{\cal G}_i$ for all $i\in S$, so $B_A\in{\cal F}$. However, $\alpha_{\cal F}({\cal P}^u)\geq l$ means that ${B_A}^c\in{\cal F}$, a contradiction.

For the other inequality, assume that the right-hand side of (\ref{eqlimpat}) is at least $l$. This means that, for every $A\in{\cal P}^u\cap{\cal C}$, the set $T:=\{i\in I:\sum_{{\cal Q}^v\in\overline{A}}\alpha_{{\cal G}_i}({\cal Q}^v)\geq l\}$ belongs to $\cal H$. If we define $B_A$ as above, we get that ${B_A}^c\in{\cal G}_i$ for $i\in T$, hence ${B_A}^c\in{\cal F}$. Thus $\sum_{{\cal Q}^v\in\overline{A}}\alpha_{{\cal F}}({\cal Q}^v)\geq l$ for all $A\in{\cal P}^u\cap{\cal C}$, so by $\cal C$-closedness of $\alpha_{\cal F}$ (Lemma \ref{cclosed}), $\alpha_{\cal F}({\cal P}^u)\geq l$ as well.\kraj

\section{$\cal C$-closedness}

We recall some more definitions from \cite{So8}.

\begin{de}
Let $(L,\leq)$ be a linear order and let $a=\langle a_m:m\in L\rangle$, $b=\langle b_m:m\in L\rangle$ be two sequences in $\mathbb{N}_\infty$. We say that $a$ dominates $b$ if, for every $l\in L$:
$$\sum_{m\geq l}a_m\geq\sum_{m\geq l}b_m.$$

For $\alpha\in{\cal A}$ and ${\cal P}\in\overline{\mathbb{P}}$, let $\alpha\upharpoonright{\cal P}=\langle {\cal P}^u:u\in{\cal E}_{\cal P}\rangle$.

Finally, for $\alpha,\beta\in{\cal A}$ we define: $\alpha\preceq\beta$ if $\beta\upharpoonright{\cal P}$ dominates $\alpha\upharpoonright{\cal P}$ for every ${\cal P}\in\overline{\mathbb{P}}$. If $\alpha\preceq\beta$ and $\beta\preceq\alpha$, we write $\alpha\approx\beta$.
\end{de}

\begin{lm}\label{maximalm}
Let $\alpha\in{\cal A}$ be $\cal C$-closed and ${\cal P}\in\overline{\mathbb{P}}$. If $\sum_{w\in{\cal E}_{\cal P}}\alpha({\cal P}^w)$ is infinite, then there is the largest $m\in{\cal E}_{\cal P}$ such that $\sum_{w\succeq_{\cal P} m}\alpha({\cal P}^w)$ is infinite, and in fact $\alpha({\cal P}^m)=\infty$.
\end{lm}

\dokaz Let $m:=\sup\{u\in{\cal E}_{\cal P}:\sum_{w\succeq_{\cal P} u}\alpha({\cal P}^w)\mbox{ is infinite}\}$, existing by Proposition \ref{supinf}. If ${\cal P}^m$ is of the first kind, then $m$ has an immediate successor $u\in{\cal E}_{\cal P}$ and $\sum_{w\succeq_{\cal P} u}\alpha({\cal P}^w)$ is finite, so $\alpha({\cal P}^m)$ must be infinite, hence its value is $\infty$.

If ${\cal P}^m$ is of the second kind, then there is a $\prec_{\cal P}$-increasing sequence $\langle u_\xi:\xi<\gamma\rangle$ in ${\cal E}_{\cal P}$ with $m=\sup_{\xi<\gamma}u_\xi$. For every $A\in{\cal P}^m\cap{\cal C}$ there is $\zeta<\gamma$ such that $\overline{A}$ contains ${\cal P}^{u_\xi}$ for all $\xi\geq\zeta$ and therefore, since $\sum_{w\succ_{\cal P} m}\alpha({\cal P}^w)$ is finite, $\sum_{{\cal Q}^v\in\overline{A}}\alpha({\cal Q}^v)\geq\sum_{u_\zeta\preceq_{\cal P}w\preceq_{\cal P}m}\alpha({\cal P}^w)$ is infinite, so by $\cal C$-closedness we have $\alpha({\cal P}^m)=\infty$.\kraj

Using this result, the proof of Lemma 3.8 from \cite{So8} is easily adjusted to obtain the following result.

\begin{pp}\label{pomdominate}
Let $\alpha,\beta\in{\cal A}$ be $\cal C$-closed and let ${\cal P}\in\overline{\mathbb{P}}$. Then $\beta\upharpoonright{\cal P}$ dominates $\alpha\upharpoonright{\cal P}$ if and only if there is a function $f_{\cal P}:\bigcup_{u\in{\cal E}_{\cal P}}(\{u\}\times\alpha({\cal P}^u))\rightarrow{\cal E}_{\cal P}$ such that $f_{\cal P}(u,i)\succeq_{\cal P}u$ for every $(u,i)\in\bigcup_{u\in{\cal E}_{\cal P}}(\{u\}\times\alpha({\cal P}^u))$ and $|f_{\cal P}^{-1}[\{v\}]|\leq\beta({\cal P}^v)$ for every $v\in{\cal E}_{\cal P}$.
\end{pp}

It was proven in \cite{So8} that ${\cal F}\widemid{\cal G}$ implies $\alpha_{\cal F}\preceq\alpha_{\cal G}$ (Corollary 3.14); on the other hand, if $\alpha_{\cal F}\preceq\beta$ and $\beta$ is $\cal U$-closed, then there is ${\cal G}\in\beta\mathbb{N}$ such that $\alpha_{\cal G}\approx\beta$ and ${\cal F}\widemid{\cal G}$ (Theorem 4.7). Question 5.2 from the same paper asked if we can strengthen the conclusion to get $\alpha_{\cal G}=\beta$. Lemma \ref{cclosed} and Example \ref{UCclosed} show that, under those assumptions, we can not: not every $\cal U$-closed pattern is equal to $\alpha_{\cal G}$ for some ultrafilter $\cal G$. We now apply Theorem \ref{limpat} to improve this result, finding the proper conditions needed for this improvement. First we tighten the relation $\preceq$ in order to make a stepping stone to a stronger result.

\begin{de}
For patterns $\alpha$ and $\beta$:

-if ${\cal P}\in\overline{\mathbb{P}}\setminus\mathbb{P}$, denote $\alpha\upharpoonright{\cal P}\leq\beta\upharpoonright{\cal P}$ if $\alpha({\cal P}^u)\leq\beta({\cal P}^u)$ for all $u\in{\cal E}_{\cal P}$, and $\alpha\upharpoonright{\cal P}<\beta\upharpoonright{\cal P}$ if also the strict inequality holds for at least one $u\in{\cal E}_{\cal P}$;

-if $p\in\mathbb{P}$, denote $\alpha\upharpoonright p\leq\beta\upharpoonright p$ if $\alpha(p^k)=1$ and $\beta(p^l)=1$ for some $k\leq l$, and $\alpha\upharpoonright p\leq\beta\upharpoonright p$ in case $k<l$;

-$\alpha\leq\beta$ if $\alpha\upharpoonright{\cal P}\leq\beta\upharpoonright{\cal P}$ for all ${\cal P}\in\overline{\mathbb{P}}$;

-$\alpha<\beta$ if $\alpha\leq\beta$ and $\alpha\upharpoonright{\cal P}<\beta\upharpoonright{\cal P}$ for some ${\cal P}\in\overline{\mathbb{P}}$.
\end{de}

\begin{lm}\label{overunderL}
Let $\beta\in{\cal A}$ be $\cal C$-closed and ${\cal F}\in\beta\mathbb{N}$.

(a) If $\alpha_{\cal F}\leq\beta$, then there is ${\cal G}\in\beta\mathbb{N}$ such that $\alpha_{\cal G}=\beta$ and ${\cal F}\widemid{\cal G}$.

(b) If $\beta\leq\alpha_{\cal F}$, then there is ${\cal G}\in\beta\mathbb{N}$ such that $\alpha_{\cal G}=\beta$ and ${\cal G}\widemid{\cal F}$.
\end{lm}

\dokaz (a) We define by recursion a $\widemid$-chain $\langle{\cal G}_\xi\rangle$ of ultrafilters such that $\alpha_{{\cal G}_\xi}\leq\beta$ and $\langle\alpha_{{\cal G}_\xi}\rangle$ is a $<$-chain of patterns, as well as a $\zvemid$-chain $\langle x_\xi\rangle$ of their generators. Begin with ${\cal G}_0:={\cal F}$ and any generator $x_0$ of $\cal F$. If at some point $\alpha_{{\cal G}_\xi}=\beta$, we are done. Assume that $\langle{\cal G}_\zeta:\zeta<\xi\rangle$ and $\langle x_\zeta:\zeta<\xi\rangle$ were already defined. 

If $\xi=\gamma+1$ and $\alpha_{{\cal G}_\gamma}<\beta$, there are two cases. First, if $\alpha(p^k)=1$ and $\beta(p^l)=1$ for some $p\in\mathbb{P}$ and $k<l$, then let $x_\xi=x_\gamma\cdot p^{l-k}$. Otherwise, there are ${\cal P}\in\overline{\mathbb{P}}\setminus\mathbb{P}$ and ${\cal P}^u\in{\cal B}$ for which $\alpha_{{\cal G}_\gamma}({\cal P}^u)<\beta({\cal P}^u)$. Take any $p^a\models{\cal P}^u$ such that $p\nzvemid x_\gamma$, and define $x_\xi=x_\gamma\cdot p^a$. If ${\cal G}_\xi:=tp(x_\xi/\mathbb{N})$, clearly we still have $\alpha_{{\cal G}_\xi}\leq\beta$ and ${\cal G}_\gamma\widemid{\cal G}_\xi$.

If $\xi$ is a limit ordinal, we choose an ultrafilter $\cal H$ on $\xi$ containing all final segments of $\xi$ and put ${\cal G}_\xi:=\lim_{\zeta\rightarrow{\cal H}}{\cal G}_\zeta$. By Proposition \ref{sublimit}, ${\cal G}_\zeta\widemid{\cal G}_\xi$ for all $\zeta<\xi$. For every ${\cal P}^u\in{\cal B}$ by ${\cal C}$-closedness of $\beta$ there is some $A\in{\cal P}^u\cap{\cal C}$ such that $\sum_{{\cal Q}^v\in\overline{A}}\beta({\cal Q}^v)=\beta({\cal P}^u)$. For every $\zeta<\xi$, since $\alpha_{{\cal G}_\zeta}\leq\beta$, we have $\sum_{{\cal Q}^v\in\overline{A}}\alpha_{{\cal G}_\zeta}({\cal Q}^v)\leq\beta({\cal P}^u)$. Now Theorem \ref{limpat} implies that $\alpha_{{\cal G}_\xi}({\cal P}^u)\leq\beta({\cal P}^u)$. By Lemma \ref{midniz}(a), we can find $x_\xi\models{\cal G}_\xi$ such that $x_\zeta\zvemid x_\xi$ for all $\zeta<\xi$.\\

(b) Similar to (a); there are two differences in the recursion step. First, in the case $\xi=\gamma+1$ instead of multiplying we divide by some $p^a\models{\cal P}^u$ such that $p^a\zvepar x_\gamma$ and $\alpha_{{\cal G}_\gamma}({\cal P}^u)>\beta({\cal P}^u)$. In the case of a limit ordinal $\xi$, we take $\cal H$ containing all initial segments of $\xi$ instead of all final segments. We use the fact that $\sum_{{\cal Q}^v\in\overline{A}}\alpha_{{\cal G}_\zeta}({\cal Q}^v)\geq\sum_{{\cal Q}^v\in\overline{A}}\beta({\cal Q}^v)\geq\beta({\cal P}^u)$ for all $A\in{\cal P}^u\cap{\cal C}$ and all $\zeta<\xi$, so by Theorem \ref{limpat} $\alpha_{{\cal G}_\xi}({\cal P}^u)\geq\beta({\cal P}^u)$.\kraj

Now we further strengthen Lemma \ref{overunderL} to obtain the affirmative answer of Question 5.2 from \cite{So8}.

\begin{te}\label{overunder}
Let $\beta\in{\cal A}$ be $\cal C$-closed and ${\cal F}\in\beta\mathbb{N}$. If $\alpha_{\cal F}\preceq\beta$, then there is ${\cal G}\in\beta\mathbb{N}$ such that $\alpha_{\cal G}=\beta$ and ${\cal F}\widemid{\cal G}$.

\end{te}

\dokaz Note first that for $p\in\mathbb{P}$ the condition $\alpha_{\cal F}\upharpoonright p\preceq\beta\upharpoonright p$ implies $\alpha_{\cal F}\upharpoonright p\leq\beta\upharpoonright p$. So we need to take care only of ${\cal P}\in\overline{\mathbb{P}}\setminus\mathbb{P}$; let $\langle{\cal P}_\gamma:\gamma<2^{\goth c}\rangle$ be an enumeration of all of them. Similarly as in Lemma \ref{overunderL}, we will define by recursion a $\widemid$-chain of ultrafilters $\langle{\cal G}_\xi:\xi<2^{\goth c}\rangle$ such that, for every $\xi<2^{\goth c}$,
\begin{eqnarray*}
& \alpha_{{\cal G}_\xi}\preceq\beta\\
& \alpha_{{\cal G}_\xi}\upharpoonright{\cal P}_\zeta\leq\beta\upharpoonright{\cal P}_\zeta & \mbox{for all }\zeta<\xi,
\end{eqnarray*}
and a $\zvemid$-chain of their generators $\langle x_\xi:\xi<2^{\goth c}\rangle$. Begin with ${\cal G}_0:={\cal F}$ and any generator $x_0$ of $\cal F$. Assume that $\langle{\cal G}_\zeta:\zeta<\xi\rangle$ and $\langle x_\zeta:\zeta<\xi\rangle$ are already defined.

If $\xi=\gamma+1$, let ${\cal P}:={\cal P}_\gamma$. If $\alpha_{{\cal G}_\gamma}\upharpoonright{\cal P}\leq\beta\upharpoonright{\cal P}$, simply let $x_\xi:=x_\gamma$ and ${\cal G}_\xi:={\cal G}_\gamma$. Otherwise, since $\alpha_{{\cal G}_\gamma}\upharpoonright{\cal P}\preceq\beta\upharpoonright{\cal P}$, by Proposition \ref{pomdominate} there is a function $f_{\cal P}:\bigcup_{u\in{\cal E}_{\cal P}}(\{u\}\times\alpha_{{\cal G}_\gamma}({\cal P}^u))\rightarrow{\cal E}_{\cal P}$ such that $f_{\cal P}(u,i)\succeq_{\cal P}u$ for every $(u,i)\in\bigcup_{u\in{\cal E}_{\cal P}}(\{u\}\times\alpha_{{\cal G}_\gamma}({\cal P}^u))$ and $|f_{\cal P}^{-1}[\{v\}]|\leq\beta({\cal P}^v)$ for every $v\in{\cal E}_{\cal P}$. We now consider three cases.



Case 1. The sum $\sum_{u\in{\cal E}_{\cal P}}\alpha_{{\cal G}_\gamma}({\cal P}^u)$ is finite. For every $u\in{\cal E}_{\cal P}$ such that $\alpha_{{\cal G}_\gamma}({\cal P}^u)>0$ we make a list of primes $p_{u,1}<p_{u,2}<\dots<p_{u,k_u}$ and $a_{u,1},a_{u,2},\dots$, $a_{u,k_u}\in\zve{\mathbb{N}}$ such that $p_{u,i}^{a_{u,i}}$ are all exact divisors of $x_\gamma$ which are generators of ${\cal P}^u$. Now choose some $b_{u,i}\geq a_{u,i}$ such that $p_{u,i}^{b_{u,i}}\models{\cal P}^{f_{\cal P}(u,i)}$, and let
$$x_\xi:=x_\gamma\cdot\prod_{\alpha_{{\cal G}_\gamma}({\cal P}^u)>0}\prod_{1\leq i\leq k_u}p_{u,i}^{b_{u,i}-a_{u,i}}.$$
If ${\cal G}_\xi:=tp(x_\xi/\mathbb{N})$, clearly $\alpha_{{\cal G}_\xi}\upharpoonright{\cal P}\leq\beta\upharpoonright{\cal P}$. Also, for all ${\cal Q}^v\in{\cal B}$ such that ${\cal Q}\neq{\cal P}$, by Lemma \ref{topology} there is a neighborhood $\overline{A}$ of ${\cal Q}^v$ not containing any of ${\cal P}^u$ or ${\cal P}^{f_{\cal P}(u,i)}$ for which the pattern value was changed, so by Theorem \ref{limpat} $\alpha_{{\cal G}_\xi}({\cal Q}^v)=\alpha_{{\cal G}_\gamma}({\cal Q}^v)$ and hence $\alpha_{{\cal G}_\xi}\upharpoonright{\cal Q}=\alpha_{{\cal G}_\gamma}\upharpoonright{\cal Q}$.

In the remaining two cases we assume that $\sum_{u\in{\cal E}_{\cal P}}\alpha_{{\cal G}_\gamma}({\cal P}^u)$ (and hence $\sum_{u\in{\cal E}_{\cal P}}\beta({\cal P}^u)$ as well) is infinite, and thus of cardinality $\infty$. By Lemma \ref{maximalm} there is a maximal $m\in{\cal E}_{\cal P}$ such that $\sum_{u\succeq_{\cal P} m}\beta({\cal P}^u)$ is infinite and, in fact, $\beta({\cal P}^m)=\infty$.

Case 2. The sum $\sum_{u\succ_{\cal P} m}\alpha_{{\cal G}_\gamma}({\cal P}^u)$ is finite. Exactly like in Case 1, for every $p_{u,i}^{a_{u,i}}$ which is a generator of some ${\cal P}^u$ for $u\succ_{\cal P} m$ we find $p_{u,i}^{b_{u,i}}\models{\cal P}^{f_{\cal P}(u,i)}$. Now we let $y_0:=x_\gamma\cdot\prod_{u\succ_{\cal P} m}\prod_{1\leq i\leq k_u}p_{u,i}^{b_{u,i}-a_{u,i}}$ and ${\cal H}_0:=tp(y_0/\mathbb{N})$. Clearly, $\alpha_{{\cal H}_0}({\cal P}^u)\leq\beta({\cal P}^u)$ for $u\succeq_{\cal P} m$. Enumerate all $p^a\models{\cal P}^w$ for $w\preceq_{\cal P} m$ which are exact divisors of $y_0$: $\langle p_\eta^{a_\eta}:\eta<\infty\rangle$. We proceed by constructing a $\widemid$-increasing sequence $\langle{\cal H}_\eta:\eta\leq\infty\rangle$ of ultrafilters and a $\zvemid$-increasing sequence $\langle y_\eta:\eta\leq\infty\rangle$ of their generators. The first step is already done. Assume that $y_\eta$ is already defined. We choose some $b_\eta\geq a_\eta$ such that $p_\eta^{b_\eta}\models{\cal P}^m$ and put $y_{\eta+1}=y_\eta\cdot p_\eta^{b_\eta-a_\eta}$ and ${\cal H}_{\eta+1}=tp(y_{\eta+1}/\mathbb{N})$. If $\eta\leq\infty$ is a limit ordinal and we have already defined $y_\zeta$ for $\zeta<\eta$, we put $[{\cal H}_\eta]:=\lim_{\zeta<\eta}[{\cal H}_\zeta]$. $y_\eta$ is obtained by applying Lemma \ref{midniz}(b). In the end we put ${\cal G}_\xi:={\cal H}_\infty$. Directly from the construction we get $\alpha_{{\cal G}_\xi}({\cal P}^w)=0$ for $w\prec_{\cal P}m$, which together with $\alpha_{{\cal G}_\xi}({\cal P}^m)=\infty=\beta({\cal P}^m)$ and $\alpha_{{\cal G}_\xi}({\cal P}^u)=\alpha_{{\cal H}_0}({\cal P}^u)\leq\beta({\cal P}^u)$ for $u\succ_{\cal P} m$ gives us $\alpha_{{\cal G}_\xi}\upharpoonright{\cal P}\leq\beta\upharpoonright{\cal P}$. Note that we have $\alpha_{{\cal G}_\xi}\upharpoonright{\cal Q}=\alpha_{{\cal G}_\gamma}\upharpoonright{\cal Q}$ for all ${\cal Q}\neq{\cal P}$; namely, for any of the other basic classes ${\cal Q}^v$ there is some $A\in{\cal U}\cap{\cal Q}^v\setminus\bigcup_{w\preceq_{\cal P} m}{\cal P}^w$, so using Theorem \ref{limpat} we easily prove by induction on $\eta$ that $\alpha_{{\cal H}_\eta}({\cal Q}^v)=\alpha_{{\cal G}_\gamma}({\cal Q}^v)$.

Case 3. The sum $s:=\sum_{u\succ_{\cal P} m}\alpha_{{\cal G}_\gamma}({\cal P}^u)$ is also infinite. Since $\sum_{u\succ_{\cal P} w}\alpha_{{\cal G}_\gamma}({\cal P}^u)$ is finite for all $w\succ_{\cal P} m$, $s$ can only equal $\aleph_0$. Let $\langle u_n:n<\omega\rangle$ be the nonincreasing sequence of elements of ${\cal E}_{\cal P}$ in which every $u\succ_{\cal P} m$ appears $\alpha_{{\cal G}_\gamma}({\cal P}^u)$-many times. We define a sequence $\langle k_n:n<\omega\rangle$ in $\zve{\mathbb{N}}$. First, let $k_0:=x_\gamma$. For each $n<\omega$, take a generator $p_n^{a_n}\models{\cal P}^{u_n}$ which is an exact factor of $x_\gamma$ (say, the $j$-th smallest such) and, as in Case 1, choose some $b_n\geq a_n$ such that $p_n^{b_n}\models {\cal P}^{f_{\cal P}(u_n,j)}$ and let $k_{n+1}:=k_n\cdot p_n^{b_n-a_n}$. Denote ${\cal K}_n=tp(k_n/\mathbb{N})$ and let $[{\cal H}_0]:=\lim_{n\rightarrow\omega}[{\cal K}_n]$. We have $\alpha_{{\cal H}_0}({\cal P}^u)\leq\beta({\cal P}^u)$ for $u\succ_{\cal P}m$ and $\alpha_{{\cal H}_0}({\cal Q}^v)=\alpha_{{\cal G}_\gamma}({\cal Q}^v)$ for ${\cal Q}\neq{\cal P}$ and all $v\in{\cal E}_{\cal Q}$. Then proceed with the construction of a sequence $\langle{\cal H}_\eta:\eta\leq\infty\rangle$ exactly as in Case 2.

Of course, in the case of a limit ordinal $\xi$ we take $[{\cal G}_\xi]:=\lim_{\zeta\rightarrow\xi}[{\cal G}_\zeta]$ and obtain $x_\xi$ from Lemma \ref{midniz}(b). Since in the $(\zeta+1)$-st step of the recursion we assured that $\alpha_{{\cal G}_{\zeta+1}}\upharpoonright{\cal P}_\zeta\leq\beta\upharpoonright{\cal P}_\zeta$ without violating either of the conditions $\alpha_{{\cal G}_{\zeta+1}}\upharpoonright{\cal Q}\leq\beta\upharpoonright{\cal Q}$ or $\alpha_{{\cal G}_{\zeta+1}}\upharpoonright{\cal Q}\preceq\beta\upharpoonright{\cal Q}$ which previously held for other prime ultrafilters, in the end we obtain ${\cal G}_{2^{\goth c}}$ for which $\alpha_{{\cal G}_{2^{\goth c}}}\leq\beta$. Hence, we apply Lemma \ref{overunderL} to ${\cal G}_{2^{\goth c}}$ and finish the proof.\kraj

\begin{co}\label{characterization}
For every $\beta\in{\cal A}$, $\beta$ is the pattern of an ultrafilter if and only if it is $\cal C$-closed.
\end{co}

\dokaz For every $\cal C$-closed pattern $\beta$ by Theorem \ref{overunder} (or even Lemma \ref{overunderL}) applied to ${\cal F}=1$ we can find an ultrafilter $\cal G$ such that $\alpha_{\cal G}=\beta$. On the other hand, by Lemma \ref{cclosed} every $\alpha_{\cal G}$ is $\cal C$-closed.\kraj

\begin{co}\label{CimpliesU}
For every $\beta\in{\cal A}$, $\cal C$-closedness implies $\cal U$-closedness.
\end{co}

\dokaz If $\beta$ is $\cal C$-closed, then there is an ultrafilter $\cal G$ such that $\alpha_{\cal G}=\beta$. But every pattern of an ultrafilter is $\cal U$-closed by \cite{So8}, Theorem 3.10.\kraj

What do our results tell us about the structure of any nonstandard $x=\prod_{n\leq z}p_n^{h(n)}\in\zve{\mathbb{N}}$? First, note that the topology on $\cal B$ introduced in Lemma \ref{topology} naturally induces a topology on $\zve{\mathbb{P}^{exp}}$ with base sets $\zve A$ for $A\in{\cal C}$. Thus we have the following.

\begin{co}\label{nonst}
Let ${\cal P}^u\in{\cal B}$. If, for every $A\in{\cal P}^u\cap{\cal C}$, $x\in\zve{\mathbb{N}}$ has infinitely many exact factors from $\zve A$, then it also has infinitely many exact factors which are generators of ${\cal P}^u$.
\end{co}

For example, if $x$ is a generator of the limit of a sequence $\langle{\cal G}_\xi\rangle$ such that, for every $A\in{\cal P}^u\cap{\cal C}$, the number of exact factors of generators of ${\cal G}_\xi$ from $\zve A$ tends to infinity, then $x$ must have infinitely many exact factors which are genrators of ${\cal P}^u$. Corollary \ref{characterization} claims that $\cal C$-closedness is, basically, the only condition restricting the structure of exact factors of a nonstandard number.

\section{Singleton classes and immediate predecessors}

In this section we turn to two problems that will show some limitations in applying the pattern method. First we consider the question: which $=_\sim$-equivalence classes are singletons? It is known (\cite{So3}, Corollary 5.10) that all classes on the first $\omega$-many levels of the $\widemid$-quasiorder are singletons. Basic classes of the first kind are clearly also examples of singletons.

\begin{te}\label{singlton}
For every ${\cal F}\in\beta\mathbb{N}$ the following conditions are equivalent:

(i) $[{\cal F}]$ is a singleton;

(ii) there are no distinct generators $x$ and $y$ of $\cal F$ such that $x\zvemid y$;

(iii) there is a $\widemid$-antichain $A\in{\cal F}$.
\end{te}

\dokaz (i)$\Rightarrow$(ii) Assume the opposite, that there are distinct generators $x$ and $y$ of $\cal F$ such that $x\zvemid y$. Let $p\in\zve{\mathbb{P}}$ be such that $p^m\zvepar x$ and $p^n\zvepar y$ for some $m<n$. Then $px$ can not be a generator of $\cal F$, since one of $x$ and $px$ belongs to $\zve\{\prod_{i\leq k}q_i^{a_i}:(\forall i\leq k)q_i\mbox{ is prime }\land 2\mid \sum_{i\leq k}a_i\}$, and the other does not. Hence, if ${\cal G}=tp(px/\mathbb{N})$, then by Proposition \ref{ekviv} ${\cal G}=_\sim{\cal F}$, a contradiction.

(ii)$\Rightarrow$(i) Assume that there is some ${\cal G}\neq{\cal F}$ such that ${\cal G}=_\sim{\cal F}$. Using Proposition \ref{ekviv}, take any generator $x$ of $\cal F$, any generator $z$ of $\cal G$ such that $x\zvemid z$ and any generator $y$ of $\cal F$ such that $z\zvemid y$. Then $x$ and $y$ contradict (ii).

(ii)$\Leftrightarrow$(iii) By Lemma 2.2.12 from \cite{L} the condition $(\exists A\in{\cal F})(\forall a,b\in A)(a\nmid b\land b\nmid a)$ is equivalent to: for all generators $x$ and $y$ of ${\cal F}$, $x\nzvemid y\land y\nzvemid x$.\kraj

Define
$${\rm supp}\alpha:=\{{\cal P}^u\in{\cal B}:\alpha({\cal P}^u)>0\}.$$

\begin{lm}\label{suppfin}
If $\alpha\in{\cal A}$ is $\cal C$-closed and $\alpha({\cal P}^u)<\infty$ for all basic classes ${\cal P}^u$, then ${\rm supp}\alpha$ is finite.
\end{lm}

\dokaz The set ${\rm supp}\alpha$ is closed in the topology introduced in Lemma \ref{topology}: if some ${\cal Q}^v\notin{\rm supp}\alpha$, then by the $\cal C$-closedness of $\alpha$ there is a neighborhood $\overline{A}$ of ${\cal Q}^v$ disjoint from ${\rm supp}\alpha$. By Lemma \ref{topology}(c), ${\rm supp}\alpha$ is therefore compact. But it is also discrete: for every ${\cal P}^u$, $\alpha({\cal P}^u)<\infty$ and $\cal C$-closedness imply that there is a neighborhood $\overline{A}$ of ${\cal P}^u$ containing no other basic classes from ${\rm supp}\alpha$. Hence ${\rm supp}\alpha$ must be finite.\kraj

\begin{te}\label{singltonpattern}
For every equivalence class $[{\cal F}]\in\beta\mathbb{N}/=_\sim$ which is not a singleton at least one of the following conditions is satisfied:

${}\hspace{5mm}$(1) there is ${\cal P}^u\in{\cal B}$ such that $\alpha_{\cal F}({\cal P}^u)=\infty$; or

${}\hspace{5mm}$(2) there is ${\cal P}^u\in{\cal B}$ of the second kind such that $\alpha_{\cal F}({\cal P}^u)>0$.
\end{te}

\dokaz Assume the opposite, that both (1) and (2) are false. By Lemma \ref{suppfin} ${\rm supp}\alpha_{\cal F}$ is finite, and it contains only basic classes of the first kind. Let $\alpha_{\cal F}=\{({\cal P}_i^{u_i},m_i):i<n\}$ and $m=\sum_{i<n}m_i$. Assume that ${\cal G}\neq{\cal F}$ is such that ${\cal G}=_\sim{\cal F}$. Since $\cal F$ contains
$$B:=\left\{\prod_{i<m}q_{i}^{a_{i}}:q_1,q_2,\dots,q_m\in\mathbb{P}\mbox{ are distinct}\right\},$$
it follows that generators of $\cal G$ must also have exactly $m$ prime factors.

We prove that actually $\alpha_{\cal F}=\alpha_{\cal G}$. For every ${\cal P}^u\in{\rm supp}\alpha_{\cal F}$ there is some $A\in{\cal P}^u\cap{\cal C}$ not contained in any other ${\cal Q}^v\in{\rm supp}\alpha_{\cal F}\cup{\rm supp}\alpha_{\cal G}$. If $\alpha_{\cal F}({\cal P}^u)=k$, then the set
$$A^{(k)}\cdot\{n\in\mathbb{N}:\neg(\exists a\in A)a\parallel n\}$$
belongs to $\cal F$, and thus also to $\cal G$, which means that $\alpha_{\cal G}({\cal P}^u)=k$ as well.

Now, let $x=p_1^{a_1}p_2^{a_2}\dots p_m^{a_m}\models{\cal F}$ and $y=p_1^{b_1}p_2^{b_2}\dots p_m^{b_m}\models{\cal G}$ be such that $x\zvemid y$. We must have $a_i\leq b_i$, but also $a_i<b_i$ for at least one $i\in\{1,2,\dots,m\}$. However, all basic classes ${\cal P}^u\in{\rm supp}\alpha_{\cal F}$ are of the first kind, so such $p_i^{a_i}$ and $p_i^{b_i}$ can not be generators of the same basic class. Thus, $\cal F$ and $\cal G$ can not have the same patterns, a contradiction.\kraj

\begin{ex}
For each of the conditions (1) and (2) from Theorem \ref{singltonpattern} it is possible to find a singleton class satisfying that condition, but not the other. First, it is obvious that, for any ${\cal P}\in\overline{\mathbb{P}}$, the class ${\cal P}^{max}$ satisfies (2) but not (1).

On the other hand, let $\langle p_n:n\in\zve{\mathbb{N}}\rangle$ be the increasing enumeration of $\zve{\mathbb{P}}$. Let $f:\mathbb{N}\rightarrow\mathbb{N}$ be any function. Now take any ultrafilter $\cal F$ containing the set $\{\prod_{i\leq n}p_i^{f(p_i)}:n\in\mathbb{N}\}$. The $=_\sim$-class $[{\cal F}]$ is not a singleton by Theorem \ref{singlton}, since any two generators of $\cal F$ are $\zvemid$-comparable. All generators of $\cal F$ are of the form $x=\prod_{i\leq z}p_i^{\zve f(p_i)}$ for some $z\in\zve{\mathbb{N}}$. Each exact factor of such an $x$ is a generator of a basic class of the first kind, and for every ${\cal P}\in\overline{\mathbb{P}}\setminus\mathbb{P}$ there is some $u\in{\cal E}_{\cal P}$ such that $\alpha_{\cal F}({\cal P}^u)=\infty$.
\end{ex}

\begin{ex}\label{exnotpattern}
The converse of the implication in Theorem \ref{singltonpattern} above is false; in fact, the condition that $[{\cal F}]$ is a singleton can not be expressed only in terms of patterns. For example, if $p\in\zve{\mathbb{P}}\setminus\mathbb{P}$, ultrafilters ${\cal F}$ and ${\cal G}$ obtained as $[{\cal F}]=\lim_{n\rightarrow\omega}[tp(2^np/\mathbb{N})]$ and ${\cal G}=tp(2^pp/\mathbb{N})$ have the same pattern $\{(2^\omega,1),({\cal P},1)\}$, where ${\cal P}=tp(p/\mathbb{N})$. 

Arbitrary generator of $\cal F$ is of the form $2^xp$. Then $2^{x-1}p$ is a generator of some ${\cal F}'\neq{\cal F}$, so since $[{\cal F}]$ is the least upper bound of $tp(2^np/\mathbb{N})$, ${\cal F}'$ and $\cal F$ must be $=_\sim$-equivalent and $[{\cal F}]$ is not a singleton. However, $\cal G$ contains the antichain $\{2^qq:q\in\mathbb{P}\}$, so by Theorem \ref{singlton} it is a singleton.
\end{ex}

The above example also shows that within the same pattern we can have different but $\widemid$-comparable $=_\sim$-classes.

Question 4.3 from \cite{So5} asked if every $=_\sim$-equivalence class $[{\cal F}]\in\beta\mathbb{N}/=_\sim$ can be represented as a limit of a $\widemid$-increasing chain of limit length. In \cite{So9} we showed that the answer is {\it no} for ultrafilters $\cal F$ belonging to basic classes ${\cal P}^u$ of the first kind. We now generalize this result.

\begin{te}\label{limitrep}
For every equivalence class $[{\cal F}]\in\beta\mathbb{N}/=_\sim$, if $[{\cal F}]$ can be represented as a limit of a strictly $\widemid$-increasing chain then $[{\cal F}]$ is not a singleton.
\end{te}

\dokaz Assume that $[{\cal F}]=\lim_{\xi\rightarrow\gamma}[{\cal F}_\xi]$ for some $\widemid$-increasing chain $\langle{\cal F}_\xi:\xi<\gamma\rangle$. If $[\cal F]$ were a singleton, by Theorem \ref{singlton} $\cal F$ would contain an antichain $A$. There must exist at least two indices $\xi,\zeta\in\gamma$ such that $A\in{\cal F}_\xi$ and $A\in{\cal F}_\zeta$. However, if $x_\xi\models{\cal F}_\xi$ and $x_\zeta\models{\cal F}_\zeta$ are such that $x_\xi\zvemid x_\zeta$, then $x_\xi,x_\zeta\in\zve A$. On the other hand $\zve A$ is an antichain by Transfer principle; a contradiction.\kraj

Another problem aiming towards better understanding of $(\beta\mathbb{N}/=_\sim,\widemid)$ is to find an equivalent condition for $[{\cal F}]$ to have an immediate predecessor. First we note that having an immediate predecessor is compatible with being a limit of a $\widemid$-increasing sequence, as the following example shows.

\begin{ex}
Let ${\cal P},{\cal Q}\in\overline{\mathbb{P}}\setminus\mathbb{P}$ be distinct and let ${\cal H}\in\beta\mathbb{N}$ contain all final segments of $\mathbb{N}$. We define ${\cal F}:=\lim_{n\rightarrow{\cal H}}{\cal P}^n$. Since multiplying by $\cal Q$ from the right is continuous, ${\cal F}\cdot{\cal Q}=\lim_{n\rightarrow{\cal H}}({\cal P}^n\cdot{\cal Q})$. This ultrafilter is both a limit of a $\widemid$-increasing sequence and has an immediate predecessor $\cal F$.
\end{ex}

One may try to prove that $[{\cal F}]$ has an immediate predecessor if and only if there is ${\cal P}^u\in{\cal B}$ of the first kind such that $0<\alpha_{\cal F}({\cal P}^u)<\infty$. However, this condition is not sufficient by the next example.

\begin{ex}
(a) The $\widemid$-maximal class $MAX$ does not have an immediate predecessor. Namely, if the $=_\sim$-class of some $\cal F$ were such a predecessor, there would exist some $p\in\mathbb{P}$ and $k\in\mathbb{N}$ such that $p^k\nwidemid{\cal F}$ (otherwise, if $\cal F$ is divisible by all natural numbers, then it belongs to $MAX$ by \cite{So4}, Lemma 4.6). But then $p\cdot{\cal F}$ is strictly $\widemid$-above $\cal F$ by \cite{So4}, Lemma 4.4. and strictly below $MAX$.

(b) The pattern $\beta=\{(p^\omega,1):p\in\mathbb{P}\}\cup\{({\cal P}^{max},\infty):{\cal P}\in\overline{\mathbb{P}}\setminus\mathbb{P}\}$ is $\cal C$-closed. By Corollary \ref{characterization} there is ${\cal G}$ such that $\alpha_{\cal G}=\beta$. Clearly, ${\cal G}\in MAX$. Let $y\models{\cal G}$ and let $q\in\zve{\mathbb{P}}$ be such that $q\nzvemid y$ and ${\cal Q}=tp(q/\mathbb{N})$. For ${\cal H}:=tp(qy/\mathbb{N})$ we have $\alpha_{\cal H}({\cal Q})=1$, but $[{\cal H}]=MAX$ does not have an immediate predecessor.
\end{ex}

Hence a reasonable sufficient condition is given by the following theorem. If ${\cal P}^u\in{\cal B}$ is of the first kind, let us denote by $u+1$ the immediate successor of $u$ in ${\cal E}_{\cal P}$.

\begin{te}
If there is ${\cal P}^u\in{\cal B}$ of the first kind such that $\alpha_{\cal F}({\cal P}^{u+1})>0$ and $\sum_{w\succeq_{\cal P}u}\alpha_{\cal F}({\cal P}^w)<\infty$, then $[{\cal F}]$ has an immediate predecessor.
\end{te}

\dokaz Let $u\in {\cal E}_{\cal P}$ satisfy the given condition and let $\sum_{w\succeq_{\cal P}u+1}\alpha_{\cal F}({\cal P}^w)=k$. Let also: $x\models{\cal F}$, $p^a\models{\cal P}^{u+1}$ such that $p^a\zvepar x$, $x=py$ and $y\models{\cal G}$. By $\cal U$-closedness of $\alpha_x$, there is a set $A\in{\cal P}^{u+1}\cap{\cal U}$ such that $A\notin{\cal Q}^v$ for all ${\cal Q}^v\in{\rm supp}\alpha_x\setminus\{{\cal P}^w:w\succeq_{\cal P}u+1\}$. Then $A^{(k)}\uparrow\in{\cal F}\cap{\cal U}\setminus{\cal G}$. Hence $[{\cal F}]\neq[{\cal G}]$, and it remains to show that $[{\cal G}]$ is an immediate predecessor of $[{\cal F}]$.

Assume the opposite, that there is some ${\cal H}\in\beta\mathbb{N}$ such that $[{\cal F}]\neq[{\cal H}]\neq[{\cal G}]$ and ${\cal G}\widemid{\cal H}\widemid{\cal F}$. By Proposition \ref{ekviv} there are some $z_1\models{\cal H}$ and $z_2\models{\cal H}$ such that $z_1\zvemid x$ and $y\zvemid z_2$. Obviously, $p^{a-1}\zvepar z_2$: $p^{a-1}\zvemid z_2$ follows from $y\zvemid z_2$, and $p^a\zvemid z_2$ would mean that $x\zvemid z_2$. In a similar way we deduce that $p^a\zvepar z_1$. For $\alpha\in{\cal A}$ and ${\cal P}^u\in{\cal B}$, denote $\alpha\upharpoonright{\cal P}^u:=\langle\alpha({\cal P}^w):w\succeq_{\cal P}u\rangle$. The condition ${\cal G}\widemid{\cal H}\widemid{\cal F}$ implies that $\alpha_{\cal F}\upharpoonright{\cal P}^{u}$ dominates $\alpha_{\cal H}\upharpoonright{\cal P}^{u}$, which in turn dominates $\alpha_{\cal G}\upharpoonright{\cal P}^{u}$. But the only difference between $\alpha_{\cal F}\upharpoonright{\cal P}^{u}$ and $\alpha_{\cal G}\upharpoonright{\cal P}^{u}$ is that $\alpha_{\cal F}({\cal P}^{u+1})=\alpha_{\cal G}({\cal P}^{u+1})+1$ and $\alpha_{\cal F}({\cal P}^{u})=\alpha_{\cal G}({\cal P}^{u})-1$, so there is no sequence between them (in sense of domination). As above, there is some $B\in{\cal P}^{u}\cap{\cal U}$ such that $B\notin{\cal Q}^v$ for all ${\cal Q}^v\in{\rm supp}\alpha_x\setminus\{{\cal P}^w:w\succeq_{\cal P}u\}$. We consider two cases.

Case 1. $\alpha_{\cal H}\upharpoonright{\cal P}^{u}=\alpha_{\cal G}\upharpoonright{\cal P}^{u}$. This means that there is some $q^b\models{\cal P}^{u+1}$ such that $q^b\zvepar x$, but $q^{b-1}\zvepar z_1$. Since $q^b\zvepar y$, we must also have $q^b\zvepar z_2$.  

Case 2. $\alpha_{\cal H}\upharpoonright{\cal P}^{u}=\alpha_{\cal F}\upharpoonright{\cal P}^{u}$. This time there is some $q^b\models{\cal P}^u$ such that $q^b\zvepar x$, $q^{b+1}\zvepar z_2$, $q^b\zvepar y$ and $q^{b}\zvepar z_1$. 

In both cases we have the same $\sum_{w\succeq_{\cal P}u}\alpha_{\cal F}({\cal P}^w)$ (hence finitely many) primes such that the exact factors from $\zve B$ of $x,y,z_1$ and $z_2$ are some powers of those primes. Suppose that $p$ is the $i$-th smallest of those primes. We can define a function $f:\mathbb{N}\rightarrow\mathbb{N}$ as follows: if $r$ is the $i$-th smallest prime such that $n$ has an exact factor from $B$ which is some power of $r$, we let $f(n)$ be such that $r^{f(n)}\parallel n$. Furthermore, since ${\cal P}^u$ is of the first kind, there is some $g:\mathbb{N}\rightarrow\mathbb{N}$ such that $\zve g(p)=a$. Now, if we define $S:=\{n\in\mathbb{N}:f(n)=g(n)\}$, then clearly $z_1\in\zve S$, but $z_2\notin\zve S$. But they are generators of the same ultrafilter, a contradiction.\kraj

However, the patterns alone can not give us an equivalent condition for the problem of immediate predecessors, as the next example shows.

\begin{ex}
Let ${\cal P},{\cal Q}\in\overline{\mathbb{P}}\setminus\mathbb{P}$ and $\beta=\{({\cal P}^{max},1),({\cal Q}^\omega,1)\}$. Let $p\models{\cal P}$ and $q\models{\cal Q}$ be such that $(p,q)$ is a tensor pair, in other words that $\zve f(q)>p$ for all $f:\mathbb{N}\rightarrow\mathbb{N}$ such that $\zve f(q)\notin\mathbb{N}$. If ${\cal F}:=tp(p^qq^p/\mathbb{N})$, then $\alpha_{\cal F}=\beta$ by \cite{So8}, Example 2.8. The $=_\sim$-equivalence class $[{\cal F}]$ is a singleton by Theorem \ref{singlton}, since the antichain $\{a^bb^a:a,b\in\mathbb{P}\land a\neq b\}$ belongs to $\cal F$. Hence $p^qq^{p-1}$ generates an ultrafilter strictly below $\cal F$, and it is clearly an immediate predecessor of $\cal F$.

Now let $s,t\in\zve{\mathbb{N}}\setminus\mathbb{N}$ be such that $p^s\models{\cal P}^{max}$, $q^t\models{\cal P}^\omega$, but $s\neq \zve f(p,q,t)$ for all $f:\mathbb{N}^3\rightarrow\mathbb{N}$. If ${\cal G}:=tp(p^sq^t/\mathbb{N})$, again we have $\alpha_{\cal G}=\beta$. Assume that $[{\cal H}]$ is an immediate predecessor of $[{\cal G}]$. Then there is some $A\in{\cal G}\cap{\cal U}\setminus{\cal H}$. Let $B$ be the set of $\mid$-minimal elements of $A$. There are some $x,y\in\zve{\mathbb{N}}$ such that $p^xq^y\in B$ and $p^{x}q^{y}\zvemid p^{s}q^{t}$. If we define a function $f$ (on an appropriate subset of $\mathbb{P}^2\times\mathbb{N}$) so that $f(a,b,k)$ is the unique $l$ such that $a^lb^k\in A$, then clearly $x=\zve f(p,q,y)$. This means that $p^{x}q^{y}\neq p^{s}q^{t}$. Thus, $tp(p^xq^y/\mathbb{N})$ is between ${\cal H}$ and ${\cal G}$: it is strictly above ${\cal H}$ because it contains $A$, and strictly below ${\cal G}$ because it contains $\{n\in\mathbb{N}:(\exists z\in A)n\mid z\}$, whose complement is in $\cal U$.
\end{ex}

What causes the limited efficiency of the pattern method for these problems? It is the fact that, for $x=\prod_{i\leq z}p_i^{h(i)}$ with more than one prime divisor, the exponents $h(i)$ may be functions of prime factors other than $p_i$ and/or of each other. For example, the difference between $p^qq^p$ and some other $p^sq^t$ in the last example can not be registered by the pattern of the ultrafilter, but it affects the properties of an ultrafilter having a singleton $=_\sim$-class or an immediate predecessor.\\

{\bf Acknowledgements.} Parts of the last section, in particular Theorem \ref{singlton} and Example \ref{exnotpattern}, were obtained in cooperation with Mauro Di Nasso, Lorenzo Luperi Baglini, Marcello Mamino, Rosario Mennuni and Mariaclara Ragosta. In a joint forthcoming paper we will consider in more detail ultrafilters whose generators have finitely many prime factors, shedding some more light on the order $(\beta\mathbb{N}/=_\sim,\widemid)$.

The author gratefully acknowledges financial support of the Science Fund of the Republic of Serbia (call IDEJE, project Set-theoretic, model-theoretic and Ramsey-theoretic phenomena in mathematical structures: similarity and diversity -- SMART, grant no.\ 7750027) and of the Ministry of Science, Technological Development and Innovation of the Republic of Serbia (grant no.\ 451-03-47/2023-01/200125).

\end{document}